\documentclass[10pt]{article}

\usepackage[dvips]{graphicx}
\usepackage{multicol}
\usepackage{enumerate}

\newtheorem{theo}{Theorem}

\newtheorem{coro}{Corollary}
\newcommand{\bb} {\mbox{\boldmath $b$}}
\newcommand{\eb} {\mbox{\boldmath $e$}}
\newcommand{\pb} {\mbox{\boldmath $p$}}
\newcommand{\rb} {\mbox{\boldmath $r$}}

\newcommand{\xb} {\mbox{\boldmath $x$}}

\newcommand{\txt}  {\textrm}
\newcommand{\hs}   {\hspace{-7mm}}

\begin{document}

\title{A numerical method for calculating the Green's function arising from electronic structure theory}

\date{February 2006}

\author
{T. Sogabe$^{1}$
T. Hoshi$^{2,3}$
S.-L. Zhang$^{1}$,
T. Fujiwara$^{2,3}$\\[2mm] 
\textit{\normalsize 
$^1$ Department of Computational Science and Engineering, 
Nagoya University,}\\[-1mm]
\textit{\normalsize Furo-cho, Chikusa-ku,Nagoya 464-8603, Japan}\\
\textit{\normalsize $^2$ Department of Applied Physics,
        The University of Tokyo,} \\[-1mm]
\textit{\normalsize Hongo, 7-3-1, Bunkyo-ku, Tokyo, 113-8656, Japan}\\
\textit{\normalsize $^3$ Core Research for Evolutional Science and Technology,}\\[-1mm]
\textit{\normalsize  Japan Science and Technology Agency (CREST-JST),}\\[-1mm]
\textit{\normalsize  4-1-8 Honcho, Kawaguchi-shi, Saitama 332-0012, Japan}
}   
 
\maketitle
\begin{abstract}
We developed a fast numerical method
for complex symmetric shifted linear systems,
which is motivated by the quantum-mechanical 
(electronic-structure) theory in nanoscale materials. 
The method is named shifted Conjugate Orthogonal Conjugate Gradient (shifted COCG) method.
The formulation is given  
and several  numerical aspects are discussed. 
\end{abstract}
\section{Introduction}
\label{sec:1}

The quantum-mechanical feature of electrons plays
a crucial role in nanoscale materials
and its mathematical foundation 
is reduced to linear-algebraic problems
with given large matrices $H$, called Hamiltonian. 
The physical properties of electrons
can be described by the Green's function $G$
that is defined as inversed matrix $G(z):= (zI-H)^{-1}$
with a complex variable $z$ whose real part corresponds to energy.
Since the standard matrix-inversion procedure
requires an impractical computational cost
in case of large matrices or in nanoscale materials,
there is a strong need for the fast solution of 
the Green's function, see the references in \cite{Takayama}.

Here we introduce a new method for calculating 
the Green's function \cite{Takayama}. Let $H$ be an $N$-by-$N$ real symmetric Hamiltonian matrix, then any element of the Green's function can be written as 
\begin{eqnarray}
   G_{ij}(z)=\eb_i^T(z I-H)^{-1}\eb_j, \label{Def:Green}
\end{eqnarray}
where $\eb_i$ denotes the $i$th unit vector, and the complex energy $z=\sigma+i\delta \in C$. Note that the $i,j$ entry of (\ref{Def:Green}) can be obtained by two steps: first, compute $(zI-H)\xb=\eb_j$, and then compute $\eb_i^T\xb$. Since an integral of $G_{ij}(z)$ with respect to $\sigma$ is needed to obtain any physical quantity \cite{Takayama}, the  numerical integration requires
$
\bigl\{(\sigma_k+i\delta )I-H\bigr\}\xb^{(k)}=\eb_j\ \txt{for}\ k=1,\ldots,M. 
$
Hence, the problem reduces to solving the following shifted linear systems with complex symmetric matrices: 
\begin{eqnarray}
(A+\sigma_k I)\xb^{(k)}=\bb \quad \txt{for}\quad k=1,\ldots, M. \label{EQ:s-systems}
\end{eqnarray}
The paper is organized as follows: in the next section, we describe the algorithm and the property of COCG for solving complex symmetric linear systems.  In \S 3, to solve (\ref{EQ:s-systems}) efficiently, a numerical method named shifted COCG is proposed and seed switching technique is introduced. In \S 4, we report some numerical experiments. Finally, we make some concluding remarks in \S 5.
                      \section{The COCG method}
Matrix $A$ is called complex symmetric if $A$ is not Hermitian but symmetric $A=A^T\neq A^H $.  To solve the linear systems, 
the COCG method \cite{COCG} has been proposed  and is known as one of the most successful Krylov subspace methods.  \\ \\
\fbox{\textbf{Algorithm 1: COCG}} \\[-5mm]
\begin{eqnarray*} 
&& \xb_0\ \txt{is an initial guess},              \\
&&\rb_{0}=\bb-A\xb_0,\ \pb_{-1}=\mathbf{0},\ \beta_{-1}=0,               \\
&& \txt{{\bf for}}\ n=0,1,\ldots\, \txt{until}\ \|\rb_n\|\le\epsilon_1\|\bb\|\
   \txt{{\bf do:}}\\                          
&&\quad \pb_{n}     = \rb_{n}+ \beta_{n-1}\pb_{n-1},                     \\
&&\quad \alpha_{n}  = \frac{\rb_n^T \rb_n}{\pb_n^T A\pb_n},              \\
&&\quad \xb_{n+1}   = \xb_n+\alpha_n\pb_n,    	                         \\
&&\quad \rb_{n+1}   = \rb_n -\alpha_n A\pb_n,                            \\
&&\quad \beta_{n}=\frac{\rb_{n+1}^T \rb_{n+1}}{\rb_n^T \rb_n}.           \\
&&\txt{{\bf end}} 
\end{eqnarray*}
Observing Algorithm 1, we see that the $n$th residual can be written as
$
  \rb_n(:=\bb-A\xb_n)=R_n(A)\rb_0,        
$
where $R_{0}(\lambda)=1,\ R_{1}(\lambda)=(1-\alpha_{0}\lambda)R_{0}(\lambda),$ 
and   
\begin{eqnarray}
R_{n}(\lambda)=\Bigl(1+\frac{\beta_{n-2}}{\alpha_{n-2}}
   \alpha_{n-1}-\alpha_{n-1}\lambda\Bigr)R_{n-1}(\lambda)    
-\frac{\beta_{n-2}}{\alpha_{n-2}}\alpha_{n-1}R_{n-2}(\lambda).	
	       			\label{3-term_3} \label{COCG_poly}
\end{eqnarray} 
It is known that if breakdown does not occur, then the $n$th residual satisfies 
\begin{eqnarray}
  \rb_n \perp \overline{K_n(A,\rb_0)},   \label{COCG_ortho}
\end{eqnarray}
which leads to conjugate orthogonality $\rb_i \perp \overline{\rb}_j$ for $i\neq j$. 
                      \section{A shifted COCG method}
In this section, we develop the COCG method for solving complex symmetric shifted linear systems.
First, let us describe the the following theorem:
\begin{theo}[{Frommer \cite[Theorem 1]{Frommer}}]
Let $W_1 \subseteq W_2\subseteq \ldots \subseteq W_k$ be a sequence of nested subspaces of $C^N$ such that $W_n$ has dimension $n$ and $W_n \cap (K_{n+1}(A,\bb))^{\perp}=\{0\},\ n=1,\ldots,k.$ Let $\rb_n:=R_n(A)\bb,\ \rb_n^{\sigma}:=R_n^{\sigma}(A+\sigma I)\bb $ be residual vectors satisfying
\begin{eqnarray}
\qquad  \rb_n, \rb_n^{\sigma}\ \perp\ W_n,\quad n=1,\ldots,k. \label{ortho}
\end{eqnarray}
Then $\rb_n$ and $\rb_n^{\sigma}$ are collinear.
\end{theo}
\begin{coro}
Let $\rb_n$ and $\rb_n^\sigma$ be the residual vectors of COCG started with $\xb_0=\xb_0^\sigma=0$. Then, $\rb_n$ and $\rb_n^\sigma$ are collinear, i.e., there exists $\pi_n^\sigma\in C$ such that $\rb_n=\pi_n^\sigma\rb_n^\sigma$.
\end{coro}
\noindent {\it Proof}.
Since it follows from (\ref{COCG_ortho}) that the COCG residuals satisfy (\ref{ortho}) with 
$
W_n=\overline{K_n(A,\bb)}=\overline{K_n(A+\sigma I,\bb)},
$
this result follows from Theorem 1.

Next, we give the formulas for computing $\rb_{n+1}^\sigma$ by using the information of $\rb_{n+1}$. It follows from the polynomial (\ref{COCG_poly}) that we have
\begin{eqnarray}
\rb_{n+1}&=&\Bigl(1+\frac{\beta_{n-1}}{\alpha_{n-1}}
 \alpha_{n}-\alpha_{n}A\Bigr)\rb_{n}   
-\frac{\beta_{n-1}}{\alpha_{n-1}}\alpha_{n}\rb_{n-1},\label{3term_rec_1} \\
 \rb_{n+1}^\sigma&=&\Bigl(1+\frac{\beta_{n-1}^\sigma}{\alpha_{n-1}^\sigma}
   \alpha_{n}^\sigma-\alpha_{n}^\sigma(A+\sigma I)\Bigr)\rb_{n}^\sigma  
-\frac{\beta_{n-1}^\sigma}{\alpha_{n-1}^\sigma}\alpha_{n}^\sigma
 \rb_{n-1}^\sigma.                                      \label{3term_rec_2}
\end{eqnarray}
Substituting the relation $\rb_{n}=\pi_n^\sigma\rb_n^\sigma$ into the previous recurrence (\ref{3term_rec_2}), we have
\begin{eqnarray}
\rb_{n+1}=\Bigl(1+\frac{\beta_{n-1}^\sigma}{\alpha_{n-1}^\sigma}
 \alpha_{n}^\sigma-\alpha_{n}^\sigma(A+\sigma I)\Bigr)
 \frac{\pi^\sigma_{n+1}}{\pi^\sigma_n}\rb_{n} 
-\frac{\beta_{n-1}^\sigma\alpha_{n}^\sigma\pi^\sigma_{n+1}}
 {\alpha_{n-1}^\sigma\pi^\sigma_{n-1}}\rb_{n-1}.        \label{3term_rec_3} 
\end{eqnarray}
To obtain the computational formula for $\rb_{n+1}^\sigma$, the three parameters $\alpha_n^\sigma$, $\beta_{n-1}^\sigma$, and $\pi_{n+1}^\sigma$ are essentially required. Hence, we give the formulas for the three parameters.
First, comparing the coefficients of $A\rb_n$ in (\ref{3term_rec_1}) and (\ref{3term_rec_3}), we find
\begin{eqnarray}
  \alpha_n^\sigma=\Bigl(\pi_n^\sigma/\pi_{n+1}^\sigma\Bigr)\alpha_n.  
  \label{SCOCG_alpha}
\end{eqnarray}
Second, comparing the coefficients of $\rb_{n-1}$ leads to
$
  \frac{\beta_{n-1}}{\alpha_{n-1}}\alpha_n=
  \frac{\beta_{n-1}^\sigma\alpha_{n}^\sigma\pi^\sigma_{n+1}}
       {\alpha_{n-1}^\sigma\pi^\sigma_{n-1}}.
$
Substituting  the result of (\ref{SCOCG_alpha}) into the previous equation, we have
\begin{eqnarray}
  \beta_{n-1}^\sigma=\Bigl(\pi_{n-1}^\sigma/\pi_{n}^\sigma\Bigr)^2\beta_{n-1}.
                                                           \label{SCOCG_beta}
\end{eqnarray}
Finally, comparing the coefficients of $\rb_n$, we find
$
 (1+\frac{\beta_{n-1}^\sigma}{\alpha_{n-1}^\sigma}
 \alpha_n^\sigma-\alpha_n^\sigma \sigma)\frac{\pi_{n+1}^\sigma}
 {\pi_n^\sigma}=1+\frac{\beta_{n-1}}{\alpha_{n-1}}\alpha_n.            
$
Substituting (\ref{SCOCG_alpha}) and (\ref{SCOCG_beta}) into the previous equation, we obtain 
\begin{eqnarray}
 \pi_{n+1}^\sigma=\Bigl(1+\frac{\beta_{n-1}}{\alpha_{n-1}}
 \alpha_n+\alpha_n \sigma\Bigr)\pi_{n}^\sigma-\frac{\beta_{n-1}}
 {\alpha_{n-1}}\alpha_n\pi_{n-1}^\sigma=R_{n+1}(-\sigma).  \label{SCOCG_pi}
\end{eqnarray}
The above formulation is based on the assumption that the seed and shifted systems are $A\xb=\bb$ and $(A+\sigma I)\xb=\bb$. Similarly, it can be readily generalized to solve $m$ shifted linear systems $(A+\sigma_i I)\xb^{(i)}=\bb$ using the seed system $(A+\sigma_s I)\xb=\bb$. The resulting algorithm is given below.\\ \\
\fbox{\textbf{Algorithm 2: Shifted COCG}} \\[-4mm]
\begin{multicols}{2}
\setlength{\columnseprule}{0.4pt}
\begin{eqnarray*} 
&&                                                \\[-11mm]
&&\hs \xb_0^{(i)}=\pb_{-1}^{(i)}=\mathbf{0},\ 
  \rb_{0}=\bb,\ \\
&&\hs \beta_{-1}=0, \ \pi^{(s, i)}_0=\pi^{(s, i)}_{-1}=\alpha_{-1}=1,    \\
&&\hs \txt{{\bf for}}\ n=0,1,\ldots\, \txt{until}\ 
         \frac{\|\rb_n\|}{\|\bb\|}\le\epsilon_1\ \txt{{\bf do:}}         \\  
&&\hs\quad \pb_{n} = \rb_{n}+ \beta_{n-1}\pb_{n-1},                      \\
&&\hs\quad \alpha_{n}  =\frac{\rb_n^T \rb_n}{\pb_n^T (A+\sigma_s I)\pb_n}, \\
&&\hs\quad  \xb_{n+1}  = \xb_n+\alpha_n\pb_n, 	                         \\
&&\hs\quad \txt{\{Begin shifted system\}}                                \\
&&\hs\quad \txt{{\bf for} $i (\neq s)=1,\ldots,m$}\ \txt{{\bf do:}}      \\
&&\hs\quad \quad \txt{\bf if $\|\rb_n^{(i)}\|> \epsilon_2\|\bb\|$}\  
              \txt{\bf then}                                             \\
&&\hs\qquad \quad  \pi^{(s,i)}_{n+1}=R_{n+1}^{(s)}(\sigma_s - \sigma_i),  
                \leftarrow (\ref{SCOCG_pi})\\
&&\hs\qquad \quad \beta_{n-1}^{(i)}=\Bigl(\pi^{(s,i)}_{n-1}/
                                  \pi^{(s,i)}_n\Bigr)^2\beta_{n-1},      \\    
&&\hs\qquad \quad \alpha_n^{(i)}=\frac{\pi^{(s,i)}_n}
                                    {\pi^{(s,i)}_{n+1}}\alpha_n,         \\
&&\hs\qquad \quad \pb_n^{(i)}=\frac{1}{\pi^{(s,i)}_n}\rb_n+
                          \beta_{n-1}^{(i)}\pb^{(i)}_{n-1},\\
&&\hs\qquad \quad \xb_{n+1}^{(i)}=\xb_n^{(i)}+\alpha_n^{(i)}\pb_n^{(i)}, \\
&&\hs\quad  \quad \txt{{\bf end if}}                                     \\
&&\hs\quad \txt{{\bf end}}                                               \\
&&\hs\quad \txt{\{End shifted system\}}                                  \\
&&\hs\quad \rb_{n+1}   = \rb_n -\alpha_n  (A+\sigma_s I)\pb_n,           \\
&&\hs\quad \beta_{n}=\frac{\rb_{n+1}^T \rb_{n+1}}{\rb_n^T \rb_n}.        \\
&&\hs\txt{{\bf end}} 
\end{eqnarray*}
\end{multicols}
\txt{}\\[-2mm]
         \textbf{Shifted COCG with seed switching technique}\\[2mm]
We can see from Algorithm 2 that if $|\pi_n^{(s,i)}|=|R_n^{(s)}(\sigma_s-\sigma_i)|\ge 1$, then $\|\rb_n^{(i)}\|\le \|\rb_n\| $. Hence, if we could find a seed system such that $|R_n^{(s)}(\sigma_s-\sigma_i)|\ge1$, then all shifted systems could be solved. However, it is extremely hard to find such system in $a$ $priori$ except some special cases discussed in \cite{Frommer}. In this section, we will avoid such problem by using the following strategy:
\begin{enumerate}[I.]
\item Choose a seed system, and then start Algorithm 2;
\item If the seed system was solved at $n$th iteration, then find the new one;
\item  Start Algorithm 2 from ($n+1$)th iteration using the new seed system. 
\end{enumerate}
In II, as one of criteria for choosing the new seed system $\tilde{s}$, we adopt  $\tilde{s} = \arg \max_{i\in I} \{\|\rb_{n}^{(i)} \| \}$, where $I$ denotes an index set of unsolved systems.
In III, we need two steps to switch the old seed system to the new one. First, compute
\begin{eqnarray*}
 \pi^{(s,\tilde{s})}_{n+1}=R_{n+1}^{(s)}(\sigma_s - \sigma_{\tilde{s}}),\quad 
\beta_{n}^{(\tilde{s})}=\Bigl(\pi^{(s,\tilde{s})}_{n}/ \pi^{(s,\tilde{s})}_{n+1}\Bigr)^2\beta_{n}
\end{eqnarray*}
for obtaining $\rb_{n+1}^{(\tilde{s})}$ and $\beta_{n}^{(\tilde{s})}\pb_{n}^{(\tilde{s})}$. Since it follows from $\rb_{n+1}^{(\tilde{s})}+\beta_{n}^{(\tilde{s})}\pb_{n}^{(\tilde{s})}$ that we have $\pb_{n+1}^{(\tilde{s})}$, we can start COCG solving the system $(A+\sigma_{\tilde{s}}I)\xb^{(\tilde{s})}=\bb$ from $(n+1)$th iteration step. Second, to solve remaining systems by using the new seed $\tilde{s}$, it requires generating $\alpha_{n+1}^{(i)},\ \beta_n^{(i)}$ from the new seed. We see that they can be readily generated by the following polynomial:
\begin{eqnarray*}
 \pi^{(\tilde{s},i)}_{n+1}=R_{n+1}^{(\tilde{s})}(\sigma_{\tilde{s}} - \sigma_i)\quad  \txt{for all}\ \ i\in I.
\end{eqnarray*}
To obtain the above polynomial, we need to compute
\begin{eqnarray*} 
 \alpha_i^{(\tilde{s})}=\Bigl(\pi^{(s,\tilde{s})}_i/\pi^{(s,\tilde{s})}_{i+1}\Bigr)\ \alpha_i, \quad
\beta_{j}^{(\tilde{s})}=\Bigl(\pi^{(s,\tilde{s})}_{j} / \pi^{(s,\tilde{s})}_{j+1}\Bigr)^2\beta_{j} 
\end{eqnarray*}
for $i=0,\ldots,n,\ j=0,\ldots,n-1$. Hence, the switching strategy requires only scalar operations, and moreover we can see that if breakdown does not occur, iterating the process from (II) to (III) enables us to keep solving the systems without losing the dimension of the Krylov subspace that has been generated until the last switching.  
                      \section{Numerical examples}
In this section, we report the results of numerical experiments. The problem originally comes from \cite{Takayama} and is written as follows: 
\begin{eqnarray*}
(\sigma_k I - H)\xb^{(k)} = \eb_1, \ \  k=1,\ldots,m,
\end{eqnarray*}
where, $\sigma_k=0.4+(k-1+i)/1000$, $H\in R^{2048\times 2048}$ is a symmetric matrix,  $\eb_1=(1,0,\ldots,0)^T$, and $m=1001$. Since $(\sigma_k I - H)$ is complex symmetric, the iterative solvers we used are COCG and shifted COCG. We can also apply shifted Bi-CGSTAB($\ell$)\cite{Frommer} and GMRES \cite{Datta} to the above problem since they can be used for general non-Hermitian shifted linear systems. However, they do not exploit the property of complex symmetric matrix. This leads to high computational costs per iteration step.

\begin{figure}
\centering
\includegraphics[height=4cm]{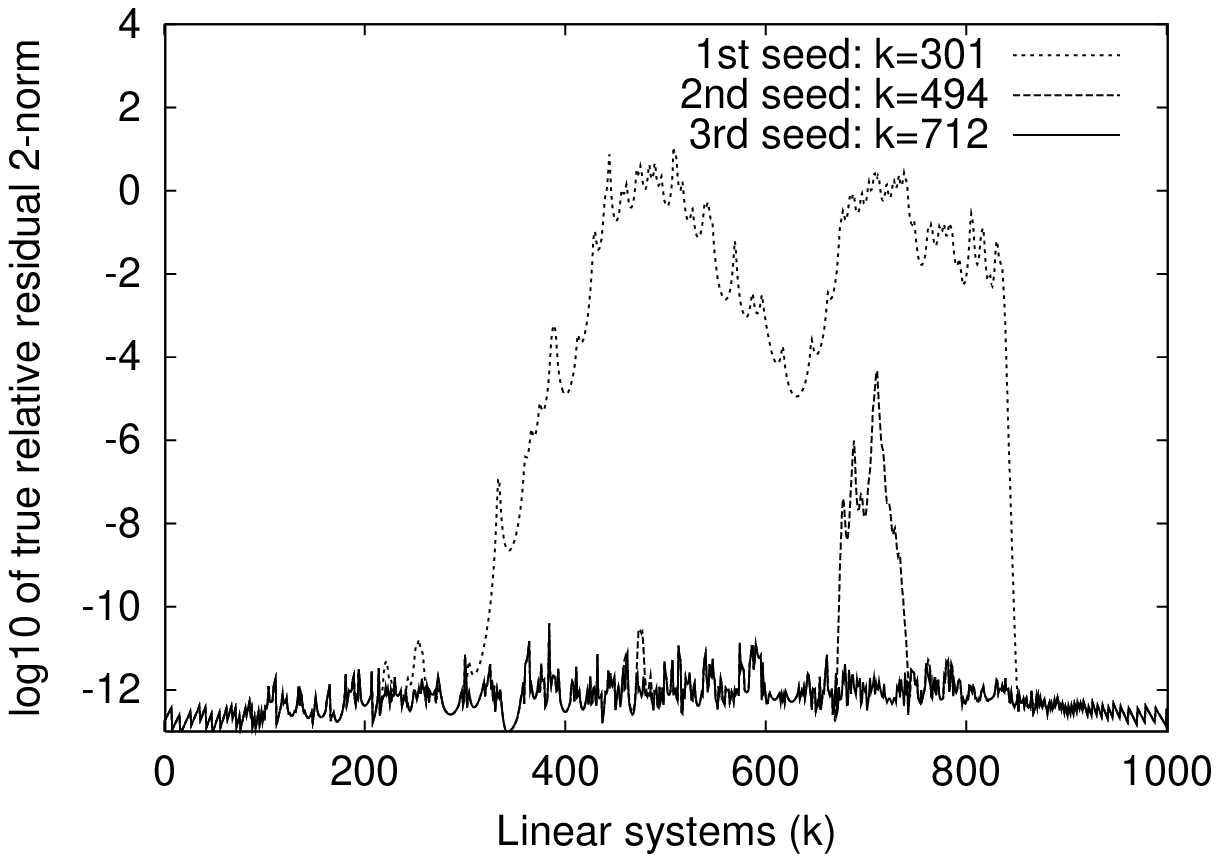}
\includegraphics[height=4cm]{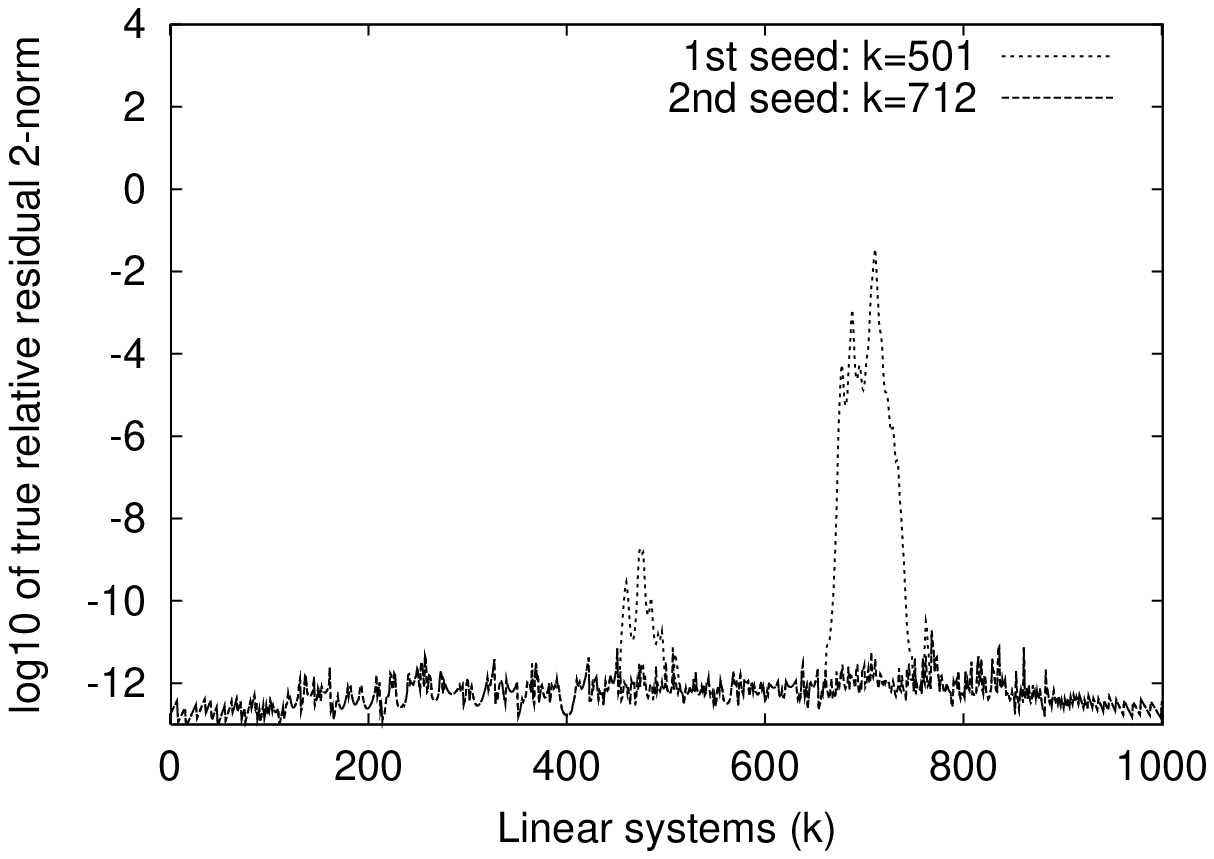}
\caption{The true relative residual 2-norm histories after each seed switching finished. The first seed is $k=301$ on the left and $k=501$ on the right.}
\label{fig:1}       
\end{figure}


\begin{table}[htb]
\begin{center}
\caption{Numerical results of COCG and shifted COCG.}
\txt{}\\[-2mm]
\label{tab:1}       
\begin{tabular}{lll}
\hline
Initial seed \ \ \  & Switching \ \ \  & Total MVs   \\
\hline
~~$k =$ 301 & ~~~~~2 & \hspace{7mm}330 \\
~~$k =$ 501 & ~~~~~1 & \hspace{7mm}328 \\
\hline
~~COCG    & ~~~~~-- & ~~124396 \\
\hline
\end{tabular}
\end{center}
\end{table}

All experiments were performed on an ALPHA work station with a 750MHz processor using double precision arithmetic. Code were written in Fortran 77 and compiled with the optimization option -$O$4.  The switching criterion is $\tilde{s} = \arg \max_{i\in I} \{\|\rb_{n}^{(i)} \| \}$. The stopping criteria are $\epsilon_1, \epsilon_2 \le 10^{-12}$.  We report two examples for $k=301$ and $k=501$ as a initial seed system.
True residual 2-norm histories are given in Fig. 1. 

In Fig. 1 on the left, 518 systems remained unsolved when the first seed system $k=301$ converged. Then, the next seed $k=494$ solved 412 more systems.  Finally, the third seed $k=712$ solved all of the remaining systems. In Fig. 1 on the right, we chose $k=501$ as an initial seed system. This choice led to 147 unsolved systems. Finally, the next seed $k=712$ solved all of the remaining systems. 

Numerical results are shown in Table 1. Total MVs denotes the total number of matrix-vector multiplications. We can see from Table 1 that shifted COCG required only about 0.27\% of Total MVs of COCG.

                      \section{Concluding remarks}
Since the original problem (\ref{Def:Green}) is a fundamental 
quantum-mechanical equation,
the present method is applicable, 
in principle,  
to various nanoscale materials, 
such as silicon, carbon, metals, polymers and so on,
if the Hamiltonian for electrons is given as an explicit matrix $H$. 

The present paper gives an interdisciplinary research between mathematics 
and physics,which shows that the computational science 
can give an important contribution to nanoscience 
through the development of general numerical algorithms,
when a fundamental equation is formulated in physics.


\end{document}